\documentclass[11pt]{article}
\usepackage[T1]{fontenc}
\usepackage[utf8]{inputenc}
\usepackage{lmodern}
\usepackage[margin=1in]{geometry}
\usepackage{amsmath,amssymb,amsthm}
\usepackage{microtype}
\usepackage{booktabs}
\usepackage{xcolor}
\usepackage[hidelinks]{hyperref}
\hypersetup{
 pdftitle={Exploring the Regular and Singular Loci of Group Representation Varieties},
 pdfauthor={S. Liriano},
 pdfsubject={Parafree groups and Borel representation varieties},
 pdfkeywords={parafree groups, representation varieties, singular loci, homology}
}
\numberwithin{equation}{section}
\newtheorem{theorem}{Theorem}
\newtheorem{proposition}[theorem]{Proposition}
\newtheorem{lemma}[theorem]{Lemma}
\newtheorem{corollary}[theorem]{Corollary}
\theoremstyle{definition}

\theoremstyle{remark}
\newtheorem{remark}[theorem]{Remark}
\newcommand{\C}{\mathbb C}
\newcommand{\Z}{\mathbb Z}
\newcommand{\R}{\mathbb R}
\newcommand{\A}{\mathbb A}

\newcommand{\Sing}{\operatorname{Sing}}
\newcommand{\Reg}{\operatorname{Reg}}
\newcommand{\SL}{\operatorname{SL}}
\newcommand{\Hom}{\operatorname{Hom}}
\newcommand{\diag}{\operatorname{diag}}
\newcommand{\rank}{\operatorname{rank}}

\makeatletter
\renewcommand{\@makefnmark}{%
  \hbox{\@textsuperscript{\normalfont\textcolor{blue}{\@thefnmark}}}}
\makeatother
\title{Exploring the Regular and Singular Loci\\
of Group Representation Varieties}
\author{S.~Liriano}
\date{}

\begin{document}
\maketitle

\begin{center}
\textit{In memory of Benjamin Fine (1948--2023).}
\end{center}

\begin{abstract}
Let \(B\) be the Borel Subgroup of \(\SL_2(\C)\), and let
\(G_m=\langle x,y,z\mid x^2y^3=z^{2m+1}\rangle\), \(m\geq2\).
The $G_m$ are parafree of rank two and deviation one. We compute the
integral homology of the regular and singular loci of
\(X_m=\Hom(G_m,B)\). Each \(X_m\) is a four-dimensional irreducible affine variety with the
homotopy type of a two-dimensional torus. Its singular locus is a disjoint
union of \(16m\) affine lines. The regular locus has Betti numbers \(b_0=1\), \(b_1=2\),
\(b_2=16m+1\), and \(b_3=b_5=16m\), with all other Betti numbers zero. Thus, although each \(X_m\) has the
homotopy type of a torus and its singular locus is a disjoint union of
contractible lines, removing the singular locus creates \(16m\) additional
independent homology classes in each of degrees \(2\), \(3\), and \(5\),
arising from the transverse link \(S^2\times S^3\). We also establish
that the regular loci are not affine.
\end{abstract}

\medskip
\noindent\textbf{Organization of the paper.}
Section~\ref{sec:introduction} states the main result and its motivation.
Section~\ref{sec:invariance} records the invariance statements and the
parafreeness of the groups under consideration.
Section~\ref{sec:equations} gives the defining equations and the
deformation retraction. Section~\ref{sec:singular} determines the singular
locus and its transverse link. Section~\ref{sec:homology} computes the
homology of the regular locus. Section~\ref{sec:consequences} gives
examples, a comparison with puncturing at a regular point, and the
nonaffineness consequence.

\section{Introduction}\label{sec:introduction}

The main aim of this paper is to demonstrate the sensitivity of the
regular and singular loci of full representation varieties as invariants
of finitely generated groups, once a suitable target algebraic group is
fixed. Our examples show that the homology of these loci distinguishes
groups whose full representation varieties have the same homotopy type
and whose affine \(B\)-conjugation quotients are isomorphic.

Representation varieties associate algebraic varieties to finitely
generated groups. Their dimension, irreducible components, and topology
provide group invariants once the target algebraic group is fixed.
The construction is classical; see Lubotzky and Magid~\cite{LM}.
The dimension sequence introduced in joint work with
Majewicz~\cite{LirianoMajewicz} records the numbers of irreducible
components in the different dimensions.

Parafree groups have the same lower central quotients as free groups of
the same rank, but need not be free. They therefore provide
a fertile ground in which to test invariants of representation varieties. The isomorphism problem for parafree groups has a
rich history.\footnote{Chandler and Magnus~\cite[p.~121]{CM}
discuss Baumslag's family
\(G(i,j)=\langle a,b,c\mid a=[c^i,a][c^j,b]\rangle\),
\(i,j\geq1\), where \([x,y]=x^{-1}y^{-1}xy\).
The question whether any two members are nonisomorphic had remained
open for several decades. Fine, Rosenberger, and Stille~\cite{FRS}
date the introduction of the family to 1962 and record the author's
1993 proof in his PhD thesis that \(G(1,1)\) and \(G(30,30)\) are not isomorphic.
The calculation involved counting points in a representation variety
over a finite field.
For the continued classification
up to isomorphism of these groups, see Chen~\cite{Chen}.}

Here the target is the Borel subgroup
\(B=\{M(\lambda,a)=
\bigl(\begin{smallmatrix}\lambda&a\\0&\lambda^{-1}\end{smallmatrix}\bigr)
:\lambda\in\C^*,\,a\in\C\}\subset\SL_2(\C)\).
For \(m\geq2\), put
\[
 G_m=\langle x,y,z\mid x^2y^3=z^{2m+1}\rangle,\qquad
 X_m=\Hom(G_m,B),\qquad
 S_m=\Sing X_m,\qquad U_m=X_m\setminus S_m.
\]
Put also \(Y=\Hom(F_2,B)=B^2\), where \(F_2\) is the free group
of rank two.\footnote{The author studied the representation varieties of the
\(G_m\) family over \(\SL_2(\C)\) in~\cite{Liriano2007}. Its representation varieties can have
an arbitrarily large number of six-dimensional irreducible components, whereas
\(\Hom(F_2,\SL_2(\C))\) has exactly one. For linear representations
of these groups, see Shalen~\cite{Shalen}.}
We regard representation varieties and their singular
loci as reduced algebraic varieties. In the sequel, homology means ordinary singular
homology with integral coefficients, in the usual complex topology.
All dimensions of complex varieties are complex dimensions.\footnote{For
the equality of algebraic and complex analytic dimension, see
Serre~\cite[\S2, no.~6, Corollary~3, p.~11]{Serre}.}

\begin{theorem}\label{thm:main}
For each \(m\geq2\), let \(N_m=16m\).
\begin{enumerate}
\item The variety \(X_m\) is an irreducible affine fourfold, and
\(X_m\simeq Y\simeq S^1\times S^1\).
\item The singular locus \(S_m\) is a disjoint union of \(N_m\) affine
lines. Along each line the local analytic model is
\[
 \C\times Q,\qquad Q=\{(s,t,u,v)\in\C^4:su+tv=0\},
\]
with the singular line corresponding to \(\C\times\{0\}\).
The link of \(Q\) at the origin is diffeomorphic to \(S^2\times S^3\).
\item The integral homology of \(U_m\) is
\begin{equation}\label{eq:main-homology}
 H_k(U_m;\Z)\cong
 \begin{cases}
 \Z,&k=0,\\
 \Z^2,&k=1,\\
 \Z^{N_m+1},&k=2,\\
 \Z^{N_m},&k=3\text{ or }5,\\
 0,&\text{otherwise}.
 \end{cases}
\end{equation}
\end{enumerate}
Consequently, if \(m\ne m'\), then \(U_m\) and \(U_{m'}\) are not
homeomorphic, and \(X_m\) and \(X_{m'}\) are not algebraically isomorphic.
The regular locus \(U_m\) is a smooth quasi-affine variety.
\end{theorem}

In this family, the resemblance to free groups extends beyond the lower
central quotients: each \(X_m\) has the homotopy type of
\(Y=\Hom(F_2,B)\). Nevertheless, the integral homology of the regular
loci distinguishes the members of the family from one another and from
the free-group case.

\begin{remark}[Affine conjugation quotient]\label{rem:borel-quotient}
For simultaneous \(B\)-conjugation and every \(m\geq2\),\footnote{Here
\(X_m/\!/B:=\operatorname{Spec}(\C[X_m]^B)\), the affine
\(B\)-conjugation quotient, as distinct from the ordinary
\(\SL_2(\C)\) character variety. In Section~\ref{sec:equations},
conjugation by \(\diag(q,q^{-1})\) fixes the diagonal entries and
sends \((a,b,c)\) to \((q^2a,q^2b,q^2c)\). Homogeneity of
\eqref{eq:X-equations} in \(a,b,c\) therefore makes
\(\C[T_m]\) the diagonal-torus invariant subring. It is also fixed
by \(B\), so \(\C[X_m]^B=\C[T_m]\cong
\C[r^{\pm1},\nu^{\pm1}]\), by Lemma~\ref{lem:torus}.}
\[
 X_m/\!/B\cong(\C^*)^2,
 \qquad
 \Sing X_m\cong\bigsqcup_{j=1}^{16m}\A^1.
\]
The affine conjugation quotient forgets the parameter \(m\), whereas
the singular locus of the full representation variety detects it.
\end{remark}

\section{Invariance and the parafree family}\label{sec:invariance}

Let \(A\) be a complex affine algebraic group and let \(G\) be a
finitely generated group. A choice of generators identifies
\(R_A(G)=\Hom(G,A)\) with the reduced algebraic set of tuples satisfying
the relations of \(G\). This is an affine variety even when an infinite
set of relations is used: its defining ideal is finitely generated
since the coordinate ring of the ambient affine variety is
Noetherian.

\begin{proposition}\label{prop:invariance}
For a fixed \(A\), the algebraic isomorphism types of \(R_A(G)\),
\(\Reg R_A(G)\), and \(\Sing R_A(G)\), and the integral homology
groups of these spaces, depend only on the isomorphism type of \(G\).
\end{proposition}

\begin{proof}
Suppose \(\varphi:G_1\xrightarrow{\sim}G_2\) is a group isomorphism.
This induces a natural pullback bijection (and in fact an isomorphism
of algebraic varieties)\footnote{Choose finite generating sets for
\(G_1\) and \(G_2\). The image under \(\varphi\) of each generator
of \(G_1\) is a word in the generators of \(G_2\). Consequently,
the coordinate maps of \(\varphi^*\) are word maps in \(A\), which
are regular because multiplication and inversion in an algebraic
group are regular. The same argument applies to
\((\varphi^{-1})^*\). Taking \(G_1=G_2=G\) and
\(\varphi=\mathrm{id}_G\) also proves independence of the chosen finite
generating set.}
\[
 \varphi^*:R_A(G_2)\longrightarrow R_A(G_1),\qquad
 \rho\longmapsto\rho\circ\varphi,
\]
with inverse \((\varphi^{-1})^*\). Thus, if \(G_1\cong G_2\), then
\(R_A(G_1)\cong R_A(G_2)\) as complex algebraic varieties.

An algebraic isomorphism induces isomorphisms of local rings and hence
preserves regular and singular points. Its restrictions to the regular
and singular loci are algebraic isomorphisms. All these maps are
homeomorphisms in the complex topology and therefore induce
isomorphisms on integral homology.
\end{proof}

The groups \(G_m\) are parafree of rank two\footnote{Recall that
\(\gamma_1(G)=G\) and \(\gamma_{r+1}(G)=[\gamma_r(G),G]\) define the
lower central series. A group \(G\) is parafree of rank \(d\) if
\(\bigcap_r\gamma_r(G)=1\) and
\(G/\gamma_r(G)\cong F_d/\gamma_r(F_d)\) for every \(r\).
Its deviation is \(d(G)-d\), where \(d(G)\) is the least number of
generators of \(G\). Abelianizing the presentation of \(G_m\) gives
\(G_m^{\mathrm{ab}}\cong
\Z^3/\langle(2,3,-(2m+1))\rangle\cong\Z^2\),
since \(\gcd(2,3,2m+1)=1\).}
and deviation one by Baumslag's results~\cite{Baumslag1967,Baumslag1968};
see also~\cite{Liriano2007}. In particular, they are not free.

\section{Defining equations and homotopy type}\label{sec:equations}

Fix \(m\geq2\) and write \(n=2m+1\). For each integer \(r\geq1\),
define the Laurent polynomial
\[
 \sigma_r(\lambda)=\sum_{j=0}^{r-1}\lambda^{r-1-2j}.
\]
Direct multiplication gives
\begin{equation}\label{eq:matrix-formulas}
 \begin{aligned}
 M(\lambda,a)M(\mu,b)
   &=M(\lambda\mu,\lambda b+a\mu^{-1}),\\
 M(\lambda,a)^r
   &=M(\lambda^r,\sigma_r(\lambda)a).
 \end{aligned}
\end{equation}
Write the images of \(x,y,z\) as \(M(\lambda,a)\),
\(M(\mu,b)\), and \(M(\nu,c)\). Then the defining equations of \(X_m\)
are
\begin{equation}\label{eq:X-equations}
 \lambda^2\mu^3=\nu^n,\qquad
 \mu^{-3}\sigma_2(\lambda)a+
 \lambda^2\sigma_3(\mu)b-\sigma_n(\nu)c=0.
\end{equation}
Let
\(T_m=\{(\lambda,\mu,\nu)\in(\C^*)^3:
\lambda^2\mu^3=\nu^n\}\).
Thus \(X_m\) is a hypersurface in \(T_m\times\C^3\), and projection
onto \(T_m\) will be denoted by \(\pi\).

\begin{lemma}\label{lem:torus}
The variety \(T_m\) is algebraically isomorphic to \((\C^*)^2\).
\end{lemma}

\begin{proof}
An isomorphism and its inverse are
\[
 (r,\nu)\longmapsto
 (r^3\nu^{-n},r^{-2}\nu^n,\nu),
 \qquad
 (\lambda,\mu,\nu)\longmapsto(\lambda\mu,\nu).
\]
Substitution in the defining equation verifies both assertions.
\end{proof}

\begin{proposition}\label{prop:irreducible}
The variety \(X_m\) is reduced and irreducible of dimension four.
\end{proposition}

\begin{proof}
By Lemma~\ref{lem:torus},
\(R=\C[T_m]\cong\C[r^{\pm1},\nu^{\pm1}]\) is a unique factorization
domain. Let \(F\in R[a,b,c]\) denote the second expression in
\eqref{eq:X-equations}. Its coefficients have no common nonunit
factor. Indeed, the common zero set in \(T_m\) of its first two
coefficients is finite: \(\sigma_2(\lambda)\) and
\(\sigma_3(\mu)\) have finitely many roots, and for each resulting
pair \((\lambda,\mu)\) there are only finitely many possible \(\nu\).
A common irreducible factor would define a divisor in the
two-dimensional torus, which cannot be contained in this finite set.

Hence \(F\) is primitive over \(R\). Over the fraction field of \(R\)
it is a nonzero polynomial of total degree one in \(a,b,c\), so it is
irreducible. Gauss's lemma\footnote{If \(A\) is a unique factorization
domain, then \(A[x_1,\ldots,x_r]\) is also a unique factorization domain.
Moreover, a nonconstant primitive polynomial in this ring is
irreducible if and only if it is irreducible over the fraction field
of \(A\).} shows that \(F\) is irreducible in \(R[a,b,c]\). In a unique factorization domain an irreducible element
is prime. Therefore \(R[a,b,c]/(F)\) is a domain and has dimension
\(2+3-1=4\).
\end{proof}

\begin{proposition}\label{prop:homotopy}
For every \(m\geq2\), \(X_m\simeq Y\simeq S^1\times S^1\).
More precisely, setting all upper-right entries equal to zero gives
a strong deformation retraction of \(X_m\) onto its diagonal
representations.
\end{proposition}

\begin{proof}
For \(0\leq\tau\leq1\), define
\(\phi_{\tau}(M(\lambda,a))=M(\lambda,(1-\tau)a)\).
The multiplication formula in \eqref{eq:matrix-formulas} shows that
each \(\phi_{\tau}:B\to B\) is a group homomorphism. Consequently,
\[
 H_{\tau}\bigl(M(\lambda,a),M(\mu,b),M(\nu,c)\bigr)
 =
 \bigl(M(\lambda,(1-\tau)a),M(\mu,(1-\tau)b),
       M(\nu,(1-\tau)c)\bigr)
\]
preserves the relation \(x^2y^3=z^n\). Equivalently, the second
equation in \eqref{eq:X-equations} is homogeneous of degree one
in \(a,b,c\), so simultaneous scaling preserves its zero set.

The map \(H:X_m\times[0,1]\to X_m\) is continuous, \(H_0\) is the
identity, \(H_1\) has image the diagonal representations, and each
\(H_{\tau}\) fixes every diagonal representation. This proves the asserted
strong deformation retraction. The diagonal representations form
\(T_m\cong(\C^*)^2\).

The same construction retracts \(Y=B^2\) onto the square of the
diagonal subgroup, again \((\C^*)^2\). Finally,
\(\zeta\mapsto\zeta|\zeta|^{-\tau}\) retracts \(\C^*\) onto \(S^1\).
\end{proof}

It follows that
\begin{equation}\label{eq:X-homology}
 H_k(X_m;\Z)\cong H_k(Y;\Z)\cong
 \begin{cases}
 \Z,&k=0\text{ or }2,\\
 \Z^2,&k=1,\\
 0,&k\geq3.
 \end{cases}
\end{equation}

\section{The singular locus and its local structure}\label{sec:singular}

The fibers of \(\pi\) are affine planes except where all three
coefficients of \(F\) vanish. Put
\[
 E_m=\{(\lambda,\mu,\nu)\in T_m:
 \sigma_2(\lambda)=\sigma_3(\mu)=\sigma_n(\nu)=0\}.
\]

\begin{lemma}\label{lem:count}
The set \(E_m\) has \(16m\) points.
\end{lemma}

\begin{proof}
For \(\zeta\ne\pm1\),
\(\sigma_r(\zeta)=\zeta^{1-r}(\zeta^{2r}-1)/(\zeta^2-1)\).
At \(\zeta=\pm1\) the value is \(r \zeta^{r-1}\), which is nonzero.
Thus the zeros of \(\sigma_r\) are precisely the \(2r\)-th roots of
unity other than \(1\) and \(-1\), and all are simple.

There are two possible \(\lambda\), each satisfying \(\lambda^2=-1\).
There are four possible \(\mu\); two satisfy \(\mu^3=1\) and two
satisfy \(\mu^3=-1\). For each pair write
\(\varepsilon=\mu^3\in\{1,-1\}\). The diagonal equation requires
\(\nu^n=-\varepsilon\). Since \(n\) is odd, exactly one of the \(n\)
roots of this equation is \(1\) or \(-1\). The other \(n-1=2m\)
roots are zeros of \(\sigma_n\). Hence
\(|E_m|=2\cdot4\cdot2m=16m\).
\end{proof}

\begin{proposition}\label{prop:local-model}
The singular locus \(S_m\) is a disjoint union of \(16m\) algebraic
affine lines, one over each point of \(E_m\). At every point of a
singular line, \(X_m\) has local analytic model \(\C\times Q\), where
\(Q=\{su+tv=0\}\subset\C^4\), and the line corresponds to
\(\C\times\{0\}\).
\end{proposition}

\begin{proof}
If a point of \(T_m\) does not lie in \(E_m\), at least one of the
partial derivatives of \(F\) with respect to \(a,b,c\) is nonzero.
Thus \(X_m\) is smooth over \(T_m\setminus E_m\).

Fix \(e=(\lambda_0,\mu_0,\nu_0)\in E_m\). The derivative with respect
to \(\nu\) of \(\lambda^2\mu^3-\nu^n\) is
\(-n\nu^{n-1}\), which is nonzero. The holomorphic implicit function
theorem therefore gives local coordinates \((\lambda,\mu)\) on
\(T_m\) near \(e\), with \(\nu\) a holomorphic function of them.
Since the zeros of \(\sigma_2\) and \(\sigma_3\) are simple,
\[
 s=\sigma_2(\lambda),\qquad t=\sigma_3(\mu)
\]
are also local coordinates, centered at \(e\).

Set \(a'=\mu^{-3}a\), \(b'=\lambda^2b\), and
\(h(s,t)=\sigma_n(\nu(s,t))\). Then \(h(0,0)=0\), and the equation
becomes \(sa'+tb'-h(s,t)c=0\). On a sufficiently small polydisc write
\[
 h(s,t)=s\alpha(s,t)+t\beta(s,t),\qquad
 \alpha(s,t)=\frac{h(s,t)-h(0,t)}{s},\quad
 \beta(s,t)=\frac{h(0,t)}{t}.
\]
Both quotients extend holomorphically across their apparent
denominators. The change of variables
\begin{equation}\label{eq:local-change}
 u=a'-\alpha(s,t)c,\qquad
 v=b'-\beta(s,t)c,\qquad w=c
\end{equation}
has inverse \(a'=u+\alpha(s,t)w\), \(b'=v+\beta(s,t)w\), \(c=w\).
The defining equation is now exactly \(su+tv=0\), independently of
\(w\).

The gradient in the four transverse variables is \((u,v,s,t)\).
It vanishes precisely at \(s=t=u=v=0\), leaving \(w\) arbitrary.
Above \(e\), this is the algebraic line
\[
 (a,b,c)=
 \bigl(\mu_0^3\alpha(0,0)w,\,
       \lambda_0^{-2}\beta(0,0)w,\,w\bigr),\qquad w\in\C.
\]
These lines are disjoint because they lie over distinct points of
\(T_m\), and there are no singular points elsewhere.
Lemma~\ref{lem:count} completes the proof.
\end{proof}

\begin{corollary}\label{cor:singular-homology}
One has \(H_0(S_m;\Z)\cong\Z^{16m}\) and
\(H_k(S_m;\Z)=0\) for \(k>0\). The reduced variety \(S_m\) is smooth;
in particular, \(\Sing(S_m)=\varnothing\).
\end{corollary}

\begin{proof}
Each component is an affine line, and hence is smooth and
contractible.
\end{proof}

\begin{lemma}\label{lem:link}
The link \(L=Q\cap S^7\), where \(S^7\) is the unit sphere in
\(\C^4\), is diffeomorphic to \(S^2\times S^3\). Consequently,
\[
 \widetilde H_j(L;\Z)\cong
 \begin{cases}
 \Z,&j=2,3,5,\\
 0,&\text{otherwise}.
 \end{cases}
\]
\end{lemma}

\begin{proof}
Make the unitary linear change of variables
\[
 \zeta_1=\frac{s+u}{\sqrt2},\quad
 \zeta_2=\frac{i(s-u)}{\sqrt2},\quad
 \zeta_3=\frac{t+v}{\sqrt2},\quad
 \zeta_4=\frac{i(t-v)}{\sqrt2}.
\]
Then \(\zeta_1^2+\zeta_2^2+\zeta_3^2+\zeta_4^2=2(su+tv)\), and the unit sphere
is preserved. Write \(\zeta=\xi+i\eta\), with \(\xi,\eta\in\R^4\).
The equations for \(L\) become
\[
 |\xi|^2=|\eta|^2=\tfrac12,\qquad \xi\cdot\eta=0.
\]
Thus \((\xi,\eta)\mapsto(\sqrt2\xi,\sqrt2\eta)\) identifies \(L\) with
the space of ordered orthonormal pairs \((a,b)\) in \(\R^4\).
This is the unit tangent bundle of \(S^3\).

Identify \(\R^4\) with the quaternions. Regard \(S^3\) as the unit
quaternions and \(S^2\) as the purely imaginary unit quaternions.
The map
\[
 S^3\times S^2\longrightarrow UT(S^3),\qquad
 (a,q)\longmapsto(a,aq)
\]
is a diffeomorphism. Its inverse is \((a,b)\mapsto(a,\bar a b)\):
orthogonality gives \(\operatorname{Re}(\bar a b)=0\), and
\(|\bar a b|=1\). This proves the claimed identification.
The homology calculation follows from the K\"unneth formula.
\end{proof}

Here \(L\) is the five-dimensional transverse link of the singular
line.\footnote{It is the link of the three-dimensional hypersurface
\(Q\), rather than the link of a point in the four-dimensional product
\(\C\times Q\). The extra factor is handled directly in the relative
homology argument below.}
The product description near a point is not the only fact needed for
that argument. We also need a neighborhood of each entire singular
line; the coordinates in \eqref{eq:local-change} provide one.

\section{Homology of the regular locus}\label{sec:homology}

\begin{lemma}\label{lem:relative}
For \(N_m=16m\), the relative homology of the pair \((X_m,U_m)\) is
\[
 H_k(X_m,U_m;\Z)\cong
 \begin{cases}
 \Z^{N_m},&k=3,4,6,\\
 0,&\text{otherwise}.
 \end{cases}
\]
\end{lemma}

\begin{proof}
Choose pairwise disjoint coordinate neighborhoods \(V_e\) of the
points \(e\in E_m\), as in Proposition~\ref{prop:local-model}.
Shrink each so that its \((s,t)\)-image is a polydisc
\(|s|,|t|<\delta_e\). The changes of the fiber variables are
invertible for all \(a,b,c\), with no restriction on their size.
Consequently the open neighborhood \(\mathcal N_e=\pi^{-1}(V_e)\) of the
entire singular line over \(e\) is analytically isomorphic, as a pair,
to
\begin{equation}\label{eq:neighborhood-pair}
 (\mathcal N_e,\mathcal N_e\setminus S_m)\cong
 \bigl(\C\times Q_{\delta_e},
       \C\times(Q_{\delta_e}\setminus\{0\})\bigr),
 \qquad
 Q_\delta=\{su+tv=0,\ |s|,|t|<\delta\}.
\end{equation}
There is no restriction on \(u,v\) in \(Q_\delta\).
Contracting the \(\C\)-factor gives a homotopy equivalence of this
pair with \((Q_{\delta_e},Q_{\delta_e}\setminus\{0\})\).

The space \(Q_\delta\) is contractible by radial contraction.
For \(0<\varepsilon<\delta\), the punctured space
\(Q_\delta\setminus\{0\}\) strongly deformation retracts onto
\(Q\cap S^7_\varepsilon\) by
\[
 p\longmapsto
 \left((1-\tau)+\tau\frac{\varepsilon}{\|p\|}\right)p,
 \qquad 0\leq\tau\leq1.
\]
This homotopy stays in \(Q_\delta\). If \(\|p\|\geq\varepsilon\),
it only decreases the absolute values of the coordinates. If
\(\|p\|<\varepsilon\), every intermediate point has norm at most
\(\varepsilon<\delta\). The homogeneous equation is preserved in
both cases. Scaling identifies \(Q\cap S^7_\varepsilon\) with \(L\).

The long exact sequence of the pair therefore gives
\[
 H_k(Q_\delta,Q_\delta\setminus\{0\};\Z)
 \cong \widetilde H_{k-1}(L;\Z)\quad(k\geq1),
 \qquad H_0(Q_\delta,Q_\delta\setminus\{0\};\Z)=0.
\]
By Lemma~\ref{lem:link}, this is \(\Z\) in degrees \(3,4,6\) and
zero otherwise.

Set \(\mathcal N=\bigcup_{e\in E_m}\mathcal N_e\). The complement \(X_m\setminus\mathcal N\)
is closed in \(X_m\) and contained in the open set \(U_m\).
Excision, followed by the finite disjoint-union decomposition, gives
\[
 H_k(X_m,U_m;\Z)
 \cong H_k(\mathcal N,\mathcal N\setminus S_m;\Z)
 \cong\bigoplus_{e\in E_m}
       H_k(\mathcal N_e,\mathcal N_e\setminus S_m;\Z).
\]
The preceding calculation and \(|E_m|=16m\) prove the result.
\end{proof}

\begin{proposition}\label{prop:regular-homology}
The groups \(H_k(U_m;\Z)\) are given by
\eqref{eq:main-homology}. In particular they are torsion-free, and
\(\chi(U_m)=-16m\).
\end{proposition}

\begin{proof}
For \(k\geq3\), equation~\eqref{eq:X-homology} gives
\(H_k(X_m;\Z)=H_{k+1}(X_m;\Z)=0\). The long exact sequence of
\((X_m,U_m)\) therefore yields
\[
 H_k(U_m;\Z)\cong H_{k+1}(X_m,U_m;\Z)\qquad(k\geq3).
\]
Lemma~\ref{lem:relative} gives
\(H_3(U_m;\Z)\cong H_5(U_m;\Z)\cong\Z^{16m}\),
\(H_4(U_m;\Z)=0\), and \(H_k(U_m;\Z)=0\) for \(k\geq6\).

In degree two, the same sequence gives
\[
 0\longrightarrow\Z^{16m}\longrightarrow H_2(U_m;\Z)
 \longrightarrow\Z\longrightarrow0.
\]
This sequence splits because \(\Z\) is free, so
\(H_2(U_m;\Z)\cong\Z^{16m+1}\).
The vanishing of the relative groups in degrees \(0,1,2\) gives
\(H_1(U_m;\Z)\cong\Z^2\) and \(H_0(U_m;\Z)\cong\Z\).
Taking the alternating sum of the ranks gives
\(1-2+(16m+1)-16m-16m=-16m\).
\end{proof}

\begin{proof}[Proof of Theorem~\ref{thm:main}]
Propositions~\ref{prop:irreducible} and~\ref{prop:homotopy} give
part~(1). Proposition~\ref{prop:local-model} and
Lemma~\ref{lem:link} give part~(2), and
Proposition~\ref{prop:regular-homology} gives part~(3).
If \(m\ne m'\), then
\(\rank H_3(U_m;\Z)=16m\ne16m'=\rank H_3(U_{m'};\Z)\), so the
regular loci cannot be homeomorphic. An algebraic isomorphism of
the whole varieties would restrict to one of their regular loci,
which is impossible. Alternatively, their singular loci have
different numbers of connected components.
The regular locus is smooth and quasi-affine because it is open in
the affine variety \(X_m\).
\end{proof}

\begin{remark}\label{rem:meskin}
We recover Meskin's result~\cite{Meskin} that the groups \(G_m\),
\(m\geq2\), are pairwise nonisomorphic. Indeed, a group isomorphism
would induce a homeomorphism of the corresponding regular loci by
Proposition~\ref{prop:invariance}, whereas \(b_3(U_m)=16m\).
\end{remark}

\section{Examples and consequences}\label{sec:consequences}

The comparison with the free group is particularly direct.
The variety \(Y=B^2\cong(\C^*)^2\times\C^2\) is a smooth affine
fourfold. Its singular locus is empty, and its homology is given by
\eqref{eq:X-homology}. Thus \(X_m\) and \(Y\) have the same homotopy
type, whereas \(U_m\) and \(Y\) do not: already their third homology
groups differ.

For the first two members of the family,
\(G_2=\langle x,y,z\mid x^2y^3=z^5\rangle\) and
\(G_3=\langle x,y,z\mid x^2y^3=z^7\rangle\), the relevant numerical
data are as follows.
\begin{center}
\setlength{\tabcolsep}{2.5pt}
\begin{tabular}{@{}lrrrrrrrr@{}}
\toprule
\(V\) & \(\#C_{\mathrm{sing}}\) & \(b_0(U)\) & \(b_1(U)\) &
 \(b_2(U)\) & \(b_3(U)\) & \(b_4(U)\) & \(b_5(U)\) & \(\chi(U)\)\\
\midrule
\(Y\)   & \(0\)   & \(1\) & \(2\) & \(1\)     & \(0\)   & \(0\) & \(0\)   & \(0\)\\
\(X_2\) & \(32\)  & \(1\) & \(2\) & \(33\)    & \(32\)  & \(0\) & \(32\)  & \(-32\)\\
\(X_3\) & \(48\)  & \(1\) & \(2\) & \(49\)    & \(48\)  & \(0\) & \(48\)  & \(-48\)\\
\(X_m\) & \(16m\) & \(1\) & \(2\) & \(16m+1\) & \(16m\) & \(0\) & \(16m\) & \(-16m\)\\
\bottomrule
\end{tabular}
\end{center}
Here \(\#C_{\mathrm{sing}}\) denotes the number of irreducible components
of \(\Sing(V)\), and \(U=V\setminus\Sing(V)\). All Betti numbers in
degrees \(k\geq6\) vanish.

\begin{remark}[Puncturing at a regular point]\label{rem:regular-puncture}
Let \(p\in U_m\) and put \(V_{m,p}=X_m\setminus\{p\}\). Then
\begin{equation}\label{eq:puncture-homology}
 H_k(V_{m,p};\Z)\cong
 \begin{cases}
 \Z,&k=0,2,7,\\
 \Z^2,&k=1,\\
 0,&\text{otherwise}.
 \end{cases}
\end{equation}
These groups are independent of \(m\) and of the chosen regular point.
For every \(q\in Y\), they are also isomorphic to
\(H_k(Y\setminus\{q\};\Z)\). The punctured variety \(V_{m,p}\)
is neither Stein nor affine.
\end{remark}

\begin{proof}
Since \(p\) is regular and \(X_m\) has complex dimension four, a
coordinate neighborhood of \(p\) is homeomorphic to an open ball in
\(\R^8\). Excision gives\footnote{See Hatcher~\cite[\S2.1,
Theorem~2.20 and the discussion of local homology on p.~126]{Hatcher}.}
\[
 H_k(X_m,V_{m,p};\Z)
 \cong H_k(\R^8,\R^8\setminus\{0\};\Z)
 \cong
 \begin{cases}
 \Z,&k=8,\\
 0,&k\ne8.
 \end{cases}
\]
By \eqref{eq:X-homology}, \(H_8(X_m;\Z)=H_7(X_m;\Z)=0\).
The long exact sequence of the pair therefore contains
\[
 0\longrightarrow\Z\xrightarrow{\ \partial\ }
 H_7(V_{m,p};\Z)\longrightarrow0.
\]
Thus \(H_7(V_{m,p};\Z)\cong\Z\), with a generator represented by a
small \(S^7\) surrounding \(p\) in the coordinate neighborhood. In
every degree \(k\ne7\), the same exact sequence shows that inclusion
induces \(H_k(V_{m,p};\Z)\cong H_k(X_m;\Z)\). This proves
\eqref{eq:puncture-homology}. The identical argument applies to the
smooth fourfold \(Y\), whose homology is also given by
\eqref{eq:X-homology}.

Finally, if \(V_{m,p}\) were Stein, Hamm's extension of the
Andreotti--Frankel Theorem to singular Stein spaces\footnote{See
Hamm~\cite{Hamm}. The singular locus \(S_m\) remains in \(V_{m,p}\),
so the singular-space version is the relevant one here.}
would give it the homotopy type of a CW complex of real dimension at
most \(\dim_{\C}V_{m,p}=4\), contradicting
\(H_7(V_{m,p};\Z)\cong\Z\). Hence \(V_{m,p}\) is not Stein.
Since complex affine varieties are Stein, it is not affine either.
\end{proof}

Removing a regular point thus creates one homology class in degree
seven, and the resulting homology does not distinguish the members
of the family or the free-group case. Removing the entire singular
locus creates \(16m\) additional independent classes in each of
degrees two, three, and five. This comparison illustrates the
sensitivity of regular-locus homology in these examples.

\begin{samepage}
\begin{corollary}\label{cor:nonaffine}
For every \(m\geq2\), the smooth fourfold \(U_m\) is not Stein and
hence is not affine.
\end{corollary}

\begin{proof}
By the Andreotti--Frankel Theorem~\cite{AF}, a Stein manifold of
complex dimension \(d\) has the homotopy type of a CW complex of real
dimension at most \(d\); see also Milnor~\cite[\S7]{Milnor} for the
smooth affine case. Since \(\dim_{\C}U_m=4\),
its fifth homology would vanish if it were Stein. But
\(H_5(U_m;\Z)\cong\Z^{16m}\ne0\).
A smooth complex affine variety is Stein, so \(U_m\) is not affine.
This does not conflict with the affineness of \(X_m\): an open
subvariety of an affine variety need not be affine, and \(X_m\)
itself satisfies Hamm's bound~\cite{Hamm}.
\end{proof}
\end{samepage}

\paragraph{Acknowledgments.}
To gain intuition, SageMath~\cite{SageMath}  through the Jupyter Notebook on a laptop
running Linux was used; no AI agent was ever employed in the process. OpenAI helped discuss aspects of  the mathematics, and in deciding if using the Borel subgroup over the reals was a good choice, something that we might revisit later. OpenAI also assisted in checking the soundness and transparency of some of my  arguments, in addition to helping
with \LaTeX{} typesetting and revising the exposition for
clarity, while keeping the presentation at an elementary level. Anthropic running in a locally hosted VM  was used to cut redundancies, and to check the references, also to reduce the number of pages in the exposition.  Google Gemini was used for locating relevant publications in the web  and downloading them to a library managed by a containerized local AI agent. An open source agent was used in converting some of my decades old notes to text. I posthumously thank Eldon Dyer for suggesting that the study of representation varieties could be approached through the study of topological invariants, something boldly underlined in my handwritten notes.

\bigskip
\noindent\small
Foro Cultural 81, Centro Hist\'orico, Guanajuato, Mexico.
Email: \href{mailto:SAL21458@yahoo.com}{SAL21458@yahoo.com}

\end{document}